\documentclass[12pt]{article}
\usepackage{amssymb,latexsym,theorem}%,txfonts}
\newcommand{\Spec}{\mathop{\mathrm{Spec}}\nolimits}

\newcommand{\ind}{\mathop{\mathrm{ind}}\nolimits}
\newcommand{\hull}{\mathop{\mathrm{hull}}\nolimits}
\newcommand{\even}{{\mathrm{even}}}
\newcommand{\hgt}{\mathop{\mathrm{ht}}\nolimits}
\newcommand{\gr}{\mathop{\mathrm{gr}}\nolimits}
\newcommand{\st}{\mathop{\mathrm{St}}\nolimits}
\newcommand{\GL}{\mathop{\mathit{GL}}\nolimits}
\newcommand{\gl}{\mathop{\mathfrak{gl}}\nolimits}

\newcommand{\SL}{\mathop{\mathit{SL}}\nolimits}

\newcommand{\Z}{\mathbb Z}
\newcommand{\Q}{\mathbb Q}
\newcommand{\Ga}{{\mathbb G_a}}
\newcommand{\Gm}{{\mathbb G_m}}
\newcommand{\C}{\mathcal C}
\newcommand{\A}{{\cal A}}
\newcommand{\qed}{\unskip\nobreak\hfill\hbox{ $\Box$}}

\newtheorem{Proposition}[subsection]{Proposition}
\newtheorem{Theorem}[subsection]{Theorem}
\newtheorem{Lemma}[subsection]{Lemma}
\newtheorem{Corollary}[subsection]{Corollary}
\theorembodyfont{\normalfont}
\newtheorem{Remark}[subsection]{Remark}
\newtheorem{Problem}[subsection]{Problem}
\newtheorem{Example}[subsection]{Example}
\newtheorem{Definition}[subsection]{Definition}
\newtheorem{Notation}[subsection]{Notation}

\begin{document}
\title{Cohomology with Grosshans graded coefficients}
\author{Wilberd van der Kallen}
\date{}
\maketitle
\sloppy
\begin{abstract}
Let the reductive
group $G$ act on the finitely generated commutative $k$-algebra
$A$. We ask if the finite generation property of the ring of invariants
$A^G$ extends
to the full cohomology ring $H^*(G,A)$. We confirm this for $G=\SL_2$
and also when the action
on $A$ is replaced by the `contracted' action on the Grosshans graded ring
$\gr A$, provided the characteristic of $k$ is large.
\end{abstract}

\section{Introduction}
Consider  a linear algebraic group or group scheme
$G$ defined over an algebraically
closed field $k$ of positive characteristic $p$.
Let $A$ be a finitely generated commutative $k$-algebra on which $G$ acts
rationally by $k$-algebra automorphisms. So $G$ acts on $\Spec(A)$.
One may then ask if the cohomology ring $H^*(G,A)$ is finitely generated
as a $k$-algebra.

We think of this as a question in invariant theory.
Our question is not what combinations of $G$ and $A$ yield a finitely generated
$H^*(G,A)$. Instead we are interested in finding
those $G$ for which every $A$
as above will give a finitely generated $H^*(G,A)$. In particular,
$G$ should be such that for all these
$A$ the ring of invariants $A^G=H^0(G,A)$
is finitely generated. This is why we must restrict attention to geometrically
reductive $G$. We believe no further restriction is needed, and our aim is
to present some evidence for this.
In characteristic zero there would be nothing to do.
Indeed, suppose $G$ is a reductive linear
algebraic group over ${\mathbb C}$.
We are concerned with
representations of the algebraic group,
the so-called
rational representations. The group is linearly reductive, meaning that all
 extensions of representations split  (even when the representations are
infinite dimensional).
So there is no higher rational cohomology and $H^*(G,A)=A^G$.
Thus our problem asks nothing new in this case. That is why we now return
to characteristic~$p$. Notice that in characteristic $p$ there is one type
of module or algebra for which the invariants, including the higher
invariants known as cohomology, behave as in characteristic zero.
They are the modules/algebras with good filtration.
They will serve as a natural tool in the sequel.

Our  proofs combine arguments and results from several earlier works.
As our question concerns `higher invariant theory', it is clear that invariant
theory will play its part. On the other hand our work
is a direct descendent of the work of Evens for the case that $G$ is a finite
group, and the work of Friedlander and Suslin for the case of finite
group schemes. We try to merge this strand
with invariant theory and emphasize that in
both cases $G$ happens to be geometrically reductive.

Thus say $G$
is geometrically reductive.
Then we know by Nagata that at least the ring of invariants $A^G=H^0(G,A)$
is finitely generated. (As explained in \cite{Borsari-Santos}, Nagata's
proof \cite{Nagata} extends to group schemes.)

By Waterhouse \cite{Waterhouse} a finite group scheme is
geometrically reductive. If $G$ is a finite
group scheme, then $A$ is a finite module over $A^G$
% \par [One shows it is integral over $A^G$, hence finite.
% For $G$ of height one, note
% that derivations fix $p$th powers. Thus one reduces to a reduced $G$,
% where \cite[4.1.2]{Springer} applies.]\par
and hence the cohomology ring $H^*(G,A)$ is indeed finitely generated by
Friedlander and Suslin \cite{Friedlander-Suslin}.
(If $G$ is finite reduced, see Evens \cite[Thm. 8.1]{Evens}.
If $G$ is finite and connected,
take $C=A^G$ in \cite[Theorem 1.5, 1.5.1]{Friedlander-Suslin}.
If $G$ is not connected, one  finishes the argument
by following \cite{Evens} as on pages 220--221 of
\cite{Friedlander-Suslin}.)

Note that if the geometrically reductive
$G$ is a subgroup of $\GL_n$, then $\GL_n/G$ is affine,
$\ind_G^{\GL_n}A$ is finitely generated, and
$H^*(G,A)=H^*(\GL_n,\ind_G^{\GL_n}A)$. (Compare \cite{Nagata},
\cite{Richardson},
\cite[Ch.\ II]{Grosshans book}, \cite[I 4.6, I 5.13]{Jantzen}.)
% \par If the groups are not  reduced, one may argue as follows.
% Take $r$ so large that the $r$-th Frobenius power $\pi=F^r:\GL_n\to \GL_n$
% maps $G$ to
% a reduced subgroup. Then $\pi(\GL_n)/\pi(G)$ is affine. And the map
% $\GL_n/G\to \pi(\GL_n)/\pi(G)$ is finite
% \cite[III,\S 3,5.5b]{Demazure-Gabriel},
% hence affine.\par
Therefore let us now assume $G=\GL_n$, or rather $G=\SL_n$ to keep it
semi-simple. (The  $\SL_n$ case suffices,
as $H^*(\GL_n,A)=H^*(\SL_n,A)^\Gm$ for a $\GL_n$-algebra $A$.)

Grosshans \cite{Grosshans contr}
has introduced a filtration on $A$. The associated graded
ring $\gr A$ is finitely generated \cite[Lemma 14]{Grosshans contr}.
There is a flat family with general fibre $A$
and special fibre $\gr A$ \cite[Theorem 13]{Grosshans contr}. Its counterpart
in characteristic zero was introduced by Popov \cite{Popov}.
Our first main result says
\begin{Theorem}\label{fingrgr}
If $n<6$ or $p>2^n$, then $H^*(\SL_n,\gr A)$ is finitely generated as a
$k$-algebra.
\end{Theorem}

\begin{Remark}
In problems \ref{integral} and \ref{divided} we discuss how
one could try to remove the
$\gr$ in the conclusion of the theorem.
\end{Remark}

Our second main result concerns the case $n=2$, where we succeed in removing
the $\gr$, using a family of universal cohomology classes which behaves as if
it is obtained by taking divided powers of the class $e_1$ of
Friedlander and Suslin
\cite{Friedlander-Suslin}.

\begin{Theorem}[Cohomological invariant theory in rank one]
Let $A$ be a finitely generated
commutative $k$-algebra on which $\SL_2$ acts rationally by algebra
automorphisms. Then $H^*(\SL_2,A)$ is finitely generated as a
$k$-algebra.
\end{Theorem}

\begin{Remark}
We know very little about the size of $H^*(G, A)$, even in simple examples.
For instance, let $p=2$ and consider the second Steinberg module $V=\st_2$
of $\SL_2$. The dimension of $V$ is four.
Its highest weight is three times the fundamental weight.
By the theorem  $H^*(\SL_2, S^*(V))$ is finitely generated,
where $S^*(V)$ denotes the symmetric algebra.
But that is all we know about the size of $H^*(\SL_2, S^*(V))$.
Note that  $S^*(V)$
does not have a good filtration
 \cite[p.71 Example]{Donkin gf}, as
$H^1(\SL_2, S^2(V))\neq0$.
\end{Remark}

\section{Recollections}
For simplicity we stay with the important case $G=\SL_n$ until
\ref{other group}.
We choose a Borel group $B^+=TU^+$ of upper triangular matrices and the
opposite Borel group $B^-$. The roots of $B^+$ are positive.
If $\lambda\in X(T)$ is dominant, then $\ind_{B^-}^G(\lambda)$ is the
dual Weyl module
$\nabla_G(\lambda)$ with highest weight $\lambda$.
In a good filtration of a $G$-module the layers are of the form
$\nabla_G(\mu)$. As in \cite{vdkallen book}
we will actually also allow a layer to be a direct sum
of any number of copies of the same $\nabla_G(\mu)$.
If $M$ is a $G$-module, and $m\geq-1$ is an integer so that
$H^{m+1}(G,\nabla_G(\mu)\otimes M)=0$ for all dominant $\mu$, then we say
as in \cite{Friedlander-Parshall}
that $M$ has \emph{good filtration dimension} at most $m$.
The case $m=0$ corresponds with $M$ having a good filtration.
We say that $M$ has good filtration dimension precisely $m$,
notation $\dim_\nabla(M)=m$,
 if $m$ is
minimal so that $M$ has good filtration dimension at most $m$.
In that case $H^{i+1}(G,\nabla_G(\mu)\otimes M)=0$ for all dominant $\mu$
and all $i\geq m$. In particular $H^{i+1}(G,M)=0$ for $i\geq m$.
If there is no finite $m$ so that  $\dim_\nabla(M)=m$, then we put
$\dim_\nabla(M)=\infty$.
\begin{Lemma}
Let $0\to M'\to M \to M''\to 0$ be exact.
Then
\begin{enumerate}
\item $\dim_\nabla(M)\leq \max(\dim_\nabla(M'),\dim_\nabla(M''))$,
\item $\dim_\nabla(M')\leq \max(\dim_\nabla(M),\dim_\nabla(M'')+1)$,
\item $\dim_\nabla(M'')\leq \max(\dim_\nabla(M),\dim_\nabla(M')-1)$,
\item $\dim_\nabla(M'\otimes M'')\leq \dim_\nabla(M')+\dim_\nabla(M'')$.
% \par In the case of type $A_n$, this is due to Wang Jian Pan
\end{enumerate}
\end{Lemma}
% Note that we do not put any restriction on the dimension of $M$.
% For us a module with $\dim_\nabla(M)$ equal to zero or $-1$ has a good
% filtration,
% regardless its dimension (compare \cite{vdkallen book}).

\subsection{Filtrations}
If $M$ is a $G$-module, and $\lambda$ is  a dominant weight,
then $M_{\leq\lambda}$ denotes the largest $G$-submodule all whose weights
$\mu$ satisfy $\mu\leq\lambda$ in the usual partial order
\cite[II 1.5]{Jantzen}.
Similarly $M_{<\lambda}$ denotes the largest $G$-submodule all whose weights
$\mu$ satisfy $\mu<\lambda$. As in \cite{vdkallen book}, we form the
$X(T)$-graded module
$$\gr_{X(T)} M=\bigoplus_{\lambda\in X(T)}M_{\leq\lambda}/M_{<\lambda}.$$
Each $M_{\leq\lambda}/M_{<\lambda}$, or $M_{\leq\lambda/<\lambda}$ for short,
 has a $B^+$-socle
$(M_{\leq\lambda/<\lambda})^U=M^U_\lambda$ of weight
$\lambda$. We view $M^U$ as a $B^-$-module through restriction (inflation)
along the
homomorphism $B^-\to T$.
Then $M_{\leq\lambda/<\lambda}$ embeds naturally in its `good filtration
hull' $\hull_\nabla(M_{\leq\lambda/<\lambda})
=\ind_{B^-}^GM^U_\lambda$.
This good filtration hull has the same $B^+$-socle
and is the injective hull
in the category $\C_\lambda$
of $G$-modules $N$ that satisfy $N=N_{\leq\lambda}$.
Compare \cite[3.1.10]{vdkallen book}.

Let us apply this in particular to our finitely generated commutative
$k$-algebra with $G$ action $A$.
We get an $X(T)$-graded algebra $\gr_{X(T)}A$. We convert it to a $\Z$-graded
algebra through an additive height function $\hgt:X(T)\to \Z$, defined by
$\hgt(\gamma)=2\sum_{\alpha>0}\langle \gamma,\widehat\alpha\rangle$,
the sum being over the positive roots.
(Our  $\hgt$ is twice the one used by Grosshans, because we prefer to get
even degrees rather than just integer degrees.)
The Grosshans graded algebra is now
$$\gr A=\bigoplus_{i\geq0}\gr_i A,$$
with $$\gr_i A=\bigoplus_{\hgt(\lambda)=i}A_{\leq\lambda/<\lambda}.$$
It embeds in a good filtration hull, which Grosshans calls $R$, and
which we call $\hull_\nabla(\gr A)$,
$$\hull_\nabla(\gr A)=\ind_{B^-}^GA^U=
\bigoplus_i\bigoplus_{\hgt(\lambda)=i}
\hull_\nabla(A_{\leq\lambda}/A_{<\lambda}).$$
Grosshans shows that $A^U$, $\gr A$, $\hull_\nabla(\gr A)$ are finitely
 generated with $\hull_\nabla(\gr A)$ finite over $\gr A$.
Mathieu did a little better in \cite{Mathieu G}.
His argument shows that in fact $\hull_\nabla(\gr A)$ is a $p$-root closure
of $\gr A$. That is,
\begin{Lemma}
For every $x\in\hull_\nabla(\gr A)$, there is an integer
$r\geq0$, so that $x^{p^r}\in \gr A$.
\end{Lemma}
\paragraph{Proof}
It suffices to take $x\in\hull_\nabla(A_{\leq\lambda}/A_{<\lambda})$
for some $\lambda$.
If $\lambda=0$, then  $\hull_\nabla(A_{\leq\lambda}/A_{<\lambda})=A^G=
\gr_0 A$. So say $\lambda>0$ and consider the subalgebra
$S=k\oplus\bigoplus_{i>0}
\hull_\nabla(A_{\leq i\lambda}/A_{<i\lambda})$ of $\hull_\nabla(\gr A)$,
with its subalgebra $S\cap\gr A$. Apply \cite[Thm 4.2.3]{vdkallen book}
to conclude that the $p$-root closure of $S\cap\gr A$ in $S$ has a good
filtration. As it contains all of $S^U$, it must be $S$ itself.\qed

\begin{Example}\label{multicone}
Consider the multicone \cite{Kempf Ramanathan}
$$k[G/U]=\ind_U^Gk=\ind_{B^+}^G\ind_U^{B^+}k=\ind_{B^+}^Gk[T]=\bigoplus_{
\lambda \textrm{ dominant}}\nabla_G(\lambda).$$
It is its own Grosshans graded ring. Recall \cite{Kempf Ramanathan} that
it is finitely generated by
the sum of the $\nabla(\varpi_i)$, where $\varpi_i$
denotes the $i$th fundamental weight. Together with the transfer principle
$A^U\cong (k[G/U]\otimes A)^G$, see \cite[Ch Two]{Grosshans book},
this gives  finite generation of $A^U$.
\end{Example}

\section{Proof of Theorem \ref{fingrgr}}
Choose $r$ so big that $x^{p^r}\in \gr A$ for all $x\in \hull_\nabla(\gr A)$.
We may view $\gr A$ as a finite $ \hull_\nabla(\gr A)^{(r)}$-module,
where the exponent $(r)$ denotes
an $r$th Frobenius twist  \cite[I 9.2]{Jantzen}.
% \par But see also the discussion in \cite{Friedlander-Suslin}, where they
% show in particular how $A$ is an $A^{(1)}$ algebra. We need the field to be
% perfect, but we do not need $A$ to be defined over the prime field.

\subsection{Key hypothesis}\label{gfhypo}
We assume that for every fundamental weight $\varpi_i$ the symmetric algebra
$S^*(\nabla(\varpi_i))$ has a good filtration.

The hypothesis in theorem \ref{fingrgr} is explained by
\begin{Lemma}
If $n<6$ or $p>2^n$, then the key hypothesis is satisfied.
\end{Lemma}
\paragraph{Proof}
We follow \cite[section 4]{Andersen-Jantzen}.
If $p>2^n$, then $p>\sum_i\dim(\nabla(\varpi_i))$, so
$\wedge^j(\nabla(\varpi_i))$ has a good filtration for all $j$ by
\cite{Andersen-Jantzen}, so
$S^*(\nabla(\varpi_i))$ has a good filtration by \cite{Andersen-Jantzen}.
If $n<6$, then by symmetry of the Dynkin diagram (contragredient duality)
it suffices to consider $\varpi_1$ and $\varpi_2$.
But $S^*(\nabla(\varpi_1))$ has a good filtration for every $n$
because $\wedge^j(\nabla(\varpi_1))$ has a good filtration for all~$j$.
And $S^*(\nabla(\varpi_2))$ has a good filtration for every $n$
by Boffi \cite{Boffi}. (Thanks to J.~Weyman and T.~Jozefiak for pointing this
out.)\qed

\begin{Example} (J.~Weyman)
If $p=2$ and $n\geq6$, there is a submodule with highest weight
$\varpi_6$
in $S^2(\nabla(\varpi_3))$.
(If $n=6$, read zero for $\varpi_6$, as we are working with $\SL_6$
rather than $\GL_6$.)
With a character computation this implies
that $S^2(\nabla(\varpi_3))$ does not have a good filtration.
The submodule in question is generated by a highest weight vector
which is a sum of ten
terms $(e_1\wedge e_{\sigma(2)}\wedge e_{\sigma(3)})
(e_{\sigma(4)}\wedge e_{\sigma(5)} \wedge e_{\sigma(6)})$.
One sums over the ten
permutations
$\sigma$ of $2,3,4,5,6$ that satisfy
$\sigma(2)<\sigma(3)$ and
$\sigma(4)<\sigma(5)<\sigma(6)$.
% \par If $p=2$ and $n\geq7$, there is a submodule with highest weight
% $\varpi_2+\varpi_6$
% in $S^3(\nabla(\varpi_3))$.
% (If $n=7$, read $\varpi_2$ for $\varpi_2+\varpi_6$.)
% With a character computation this implies
% that $S^3(\nabla(\varpi_3))$ does not have a good filtration.
% The submodule in question is generated by a highest weight vector
% which is a sum of thirty
% terms $(e_1\wedge e_{2}\wedge e_{\sigma(3)})
% (e_1\wedge e_{\sigma(4)}\wedge e_{\sigma(5)})
% (e_2\wedge e_{\sigma(6)}\wedge e_{\sigma(7)})$.
% One sums over the thirty
% permutations
% $\sigma$ of $3,4,5,6,7$ that satisfy
% $\sigma(4)<\sigma(5)$ and
% $\sigma(6)<\sigma(7)$.
\end{Example}

\paragraph{}We want to view $S^*(\nabla(\varpi_i))$ as a
\emph{graded polynomial $G$-algebra with good filtration}. Let us collect
the properties that we will use in a rather artificial definition.
\begin{Definition}
Let $D$ be a diagonalizable group scheme \cite[I 2.5]{Jantzen}.
We say that $P$ is a graded polynomial $G\times D$-algebra with good
filtration,
if the following holds. First of all $P$ is a polynomial algebra over $k$
in finitely many variables. Secondly, these variables are homogeneous of
non-negative integer degree, thus making $P$ into a graded $k$-algebra.
There is also given an action of $G\times D$ on $P$ by algebra automorphisms
that are compatible
with the grading, with $G\times D$ acting trivially on the degree zero part
$P_0$ of $P$. This $P_0$ is thus the polynomial algebra generated by the
variables of degree zero. The variables of positive degree generate their
own polynomial algebra $P^c$, which we also assume to be $G\times D$
invariant. Thus $P=P_0\otimes_kP^c$.
And finally, $P^c$ is an algebra
with a good filtration
for the action of $G$. Then of course, so is~$P$.
\end{Definition}

\begin{Example}Our key hypothesis makes that one gets a
\emph{graded polynomial $G\times D$-algebra with good filtration}
by taking for $D$ any subgroup scheme of $T$,
with $T$ acting on $P=S^*(\nabla(\varpi_i))$ through its natural
$X(T)$-grading: On $S^j(\nabla(\varpi_i))$ we make $T$ act with weight
$j\varpi_i$.
\end{Example}

If $P_i$ are graded polynomial $G\times D_i$-algebras with good filtration
for $i=1$, $2$,
then $P_1\otimes P_2$ is a graded polynomial
$G\times (D_1\times D_2)$-algebra with good filtration.

\begin{Definition}
If $P$ is a graded polynomial $G\times D$-algebra with good filtration,
then by a \emph{finite graded $P$-module} $M$ we mean a finitely generated
$\Q$-graded module for
the graded polynomial ring, together with a $G\times D$ action on $M$
which is
compatible with the grading and with
the action on $P$.
It is not required that the action is trivial on the degree zero part of $M$.
We call $M$ \emph{free} if it is free as a module over the polynomial ring.
We call $M$ \emph{extended}
if there is a finitely generated $\Q$-graded $P_0$-module $V$ with
$G\times D$
action,
so that $M=V\otimes_{P_0} P$ as graded $P$-modules.
 \end{Definition}

\begin{Lemma}\label{fresolve}
Let $P$ be a graded polynomial $G\times D$-algebra with good filtration,
and let $M$ be a finite graded $P$-module.
\begin{enumerate}
\item
There is a finite free resolution
$$0\to F_s\to F_{s-1}\to \cdots \to F_0\to M.$$
\item Every finite free $P$-module has a finite filtration whose quotients are
extended.
\item $\dim_\nabla(M)<\infty$.
\item $H^i(G,M)$ is a finite $P^G$-module for every $i\geq0$.
% \par So we are saying that at least in this case $H^*(G,M)$ is a finite
% $H^*(G,P)$-module.
\end{enumerate}
\end{Lemma}
\paragraph{Proof}  Take a finite dimensional graded
$G\times D$-submodule $V$ of $M$ that generates $M$ as a $P$-module.
Then $F_0=V\otimes_{k} P$ is free and it maps onto $M$.
As $M$ has finite projective dimension \cite[18C]{Matsumura},
we may repeat until the syzygy
is projective. But a projective module
over a polynomial ring is free by Quillen and Suslin \cite{Lam}.

Now consider a free $P$-module $F$. Let $P_+$ be the ideal generated by the
$P_i$ with $i>0$. If $F$ is free of rank $t$, then $F/P_+F$ is free of rank
$t$ over $P_0$. Choose homogeneous elements $e_1,\ldots,e_t$ in $F$ so that
their
classes form a basis of $F/P_+F$. Then the $e_i$ generate $F$
 (cf.~\cite[6.13 Lemma 5]{Jacobson}), so they also
form a basis of $F$. Let $F_m$ be the component of lowest degree in $F$,
assuming $F\neq0$. The $e_i$ that lie in $F_m$ form a basis of $F_m$ over
$P_0$. So $F_m$ generates a free, extended submodule $PF_m$ of $F$ and
the quotient $F/PF_m$ is free of lower rank.

To show that $\dim_\nabla(M)<\infty$, it suffices to
consider the case of an extended module $M=V\otimes_{P_0} P$.
As $V\otimes_{P_0} P=V\otimes_kP^c$, it suffices
to check that $V$ has finite good filtration dimension.
But  $V$ is a finitely generated $P_0$-module
and thus has only finitely many weights. Therefore the argument used
in \cite{Friedlander-Parshall} to show that finite dimensional $G$
modules have finite good filtration dimension, applies to $V$.

Finally let $M$ be  any finite graded $P$-module again.
As $G$ is reductive by Haboush \cite[10.7]{Jantzen}, it is well known that
$H^0(G,M)$ is a finite $P^G$-module,
because $(S_P^*(M))^G$ is finitely generated. So we argue by dimension shift.
We claim that
for $s$ sufficiently large the module $M\otimes \st_s\otimes \st_s$
has a good filtration.
It suffices to check this for the extended case $M=V\otimes_{P_0} P$
and then one can use again that $V$ has only finitely many weights,
so that one may choose $s$ so large that all weights of
$V\otimes k_{-(p^s-1)\rho}$ are anti-dominant. Then
$V\otimes \st_s=\ind_{B^+}^G(V\otimes k_{-(p^s-1)\rho})$ has a good filtration
and
so does $M\otimes \st_s\otimes \st_s$.
Then $H^i(G,M)$ is the cokernel of $H^{i-1}(G,M\otimes \st_s\otimes \st_s)
\to H^{i-1}(G,M\otimes \st_s\otimes \st_s/M) $ for $i\geq1$.\qed

\paragraph{} Recall that we choose $r$ so big that $x^{p^r}\in \gr A$ for all
$x\in \hull_\nabla(\gr A)$. Let $G_r$ denote the $r$-th Frobenius kernel.
We will need multiplicative structure on a Hochschild--Serre spectral sequence.
See for instance
\cite{Benson II}, translating from modules over a Hopf algebra to
comodules over a Hopf algebra.

\begin{Proposition}\label{stops}
Assume the key hypothesis \ref{gfhypo}.
With $r$ as indicated,  the spectral sequence
$$E_2^{ij}=H^i(G/G_r,H^j(G_r,\gr A))\Rightarrow H^{i+j}(G,\gr A)$$ stops,
i.e. $E_s=E_\infty$ for some $s<\infty$.
\\
In fact, $H^*(G_r,\gr A)^{(-r)}$ has finite good filtration dimension.
\\
Moreover, $E_2^{**}$ is a finite module over the even part
of~$E_2^{0*}$.
\end{Proposition}
\paragraph{Proof}
By \cite[Th 1.5, 1.5.1]{Friedlander-Suslin} the
ring $H^*(G_r,\gr A)^{(-r)}$ is a finite module
over the algebra $$\bigotimes_{a=1}^rS^*((\gl_n)^\#(2p^{a-1}))\otimes
\hull_\nabla(\gr A).$$
Here $(\gl_n)^\#(2p^{a-1})$ denotes the dual of $\gl_n$
placed in degree $2p^{a-1}$.
As $G=\SL_n$,  it follows from
\cite[4.3]{Andersen-Jantzen} that
the algebra $\bigotimes_{a=1}^rS^*((\gl_n)^\#(2p^{a-1}))$ is
a graded polynomial $G\times \{1\}$-algebra with good filtration.
So we now try to replace $\hull_\nabla(\gr A)$ by a similarly
nice algebra. That is, we seek a graded polynomial $G\times D$-algebra $P$
with good filtration and a surjection $P^D\to \hull_\nabla(\gr A)$
of graded $G$-algebras. This is where the key hypothesis comes in.
Let $v$ be a weight vector in $\hull_\nabla(\gr A)^U=A^U$.
If $v$ has weight zero, we may map $x$ to $v$, where $x$ is the variable
in a polynomial ring $k[x]$ with trivial grading and trivial $G\times T$
action.
If $v$ has weight $\lambda\neq 0$, take for $D$ the scheme
theoretic kernel of $ \lambda$ and observe
that by our hypothesis the $X(T)$-graded algebra
$$P=\bigotimes_{i=1}^{n-1}S^*(\nabla(\varpi_i))$$
is a graded polynomial $G\times D$-algebra
with good filtration, if we give $\nabla(\varpi_i)$ degree $\hgt(\varpi_i)$.
We have $P^D=\bigoplus_jP_{j\lambda}$, where
$P_{j\lambda}$ denotes the summand
$\bigotimes_{i=1}^{n-1}S^{m_i}(\nabla(\varpi_i))$ with $\sum_im_i\varpi_i=
j\lambda$. Choose a weight vector $x$ of weight $\lambda$
in $P_{\lambda}$
 (for the $G$ action).
The map from the polynomial ring $k[x]$ to
$\hull_\nabla(\gr A)^U$ which sends $x$ to $v$ extends uniquely
to a $G$ equivariant
algebra
map $P^D\to \hull_\nabla(\gr A)$ because
$P_{j\lambda}=(P_{j\lambda})_{\leq j\lambda}$.
(The first subscript in the right hand side
refers to a $T$ action associated
with the $X(T)$-grading, the second to the
$G$ action on  $P$.) Compare \cite[4.2.4]{vdkallen book}.
As $A^U$ is finitely generated, we may combine finitely many such maps
$P(i)^{D(i)}\to \hull_\nabla(\gr A)$ into one surjective map
$P^D\to \hull_\nabla(\gr A)$, with $D$ the product of the $D(i)$ and
$P$ the tensor product of the $P(i)$.
If we let the linearly reductive
$D$ act on $P\otimes_{P^D}H^*(G_r,\gr A)^{(-r)}$ through its action on $P$,
then $H^*(G_r,\gr A)^{(-r)}$ is just the direct summand, as a
$G$-module, consisting
of the $D$-invariants. So in order to show that $H^*(G_r,\gr A)^{(-r)}$
has finite good filtration dimension, it suffices to show that
$P\otimes_{P^D}H^*(G_r,\gr A)^{(-r)}$ has finite good filtration dimension.
We now view $P\otimes\bigotimes_{a=1}^rS^*((\gl_n)^\#(2p^{a-1}))$ as a
graded polynomial $G\times D$-algebra with good filtration.
(Collect the bigrading into a single total grading.)
The algebra $P\otimes_{P^D}H^*(G_r,\gr A)^{(-r)}$
is a finite graded
$P\otimes\bigotimes_{a=1}^rS^*((\gl_n)^\#(2p^{a-1}))$-module.
(We need a $\Q$-grading on  $P\otimes_{P^D}H^*(G_r,\gr A)^{(-r)}$
because of the twist.)
By Lemma \ref{fresolve} such a $P$-module
has finite good filtration dimension.
Therefore there are only finitely many $i$ with $E_2^{i*}\neq 0$
and the spectral sequence stops.
By the same lemma $H^i(G,P\otimes_{P^D}H^*(G_r,\gr A)^{(-r)})$
is finite over
$H^0(G,P\otimes\bigotimes_{a=1}^rS^*((\gl_n)^\#(2p^{a-1})))$.
Taking $D$-invariants again, we see that $E_2^{**}$ is a finite module over
$H^0(G/G_r,(P^D\otimes\bigotimes_{a=1}^rS^*((\gl_n)^\#(2p^{a-1})))^{(r)})$.
But this ring acts by way of the the even part
of~$E_2^{0*}$.
\qed

\subsection{End of proof of theorem \ref{fingrgr}}
Look some more at the spectral sequence of Proposition \ref{stops}.
We argue partly
as in Evens' proof of his finite generation theorem
\cite[4.2]{Benson II}, \cite{Evens}.
As we have already observed, it follows from \cite{Friedlander-Suslin} that
  $H^*(G_r,\gr A)$
is noetherian
 over its finitely
generated even degree part. By the proposition
$E_2$ is noetherian
 over the even degree part of $E_2^{0*}$, which is  finitely
generated because $G/G_r$ is reductive.
Or, recall from the proof of the proposition that
 the $E_2$ term is finite over the
finitely generated $k$-algebra
$H^0(G/G_r,(P^D\otimes\bigotimes_{a=1}^rS^*((\gl_n)^\#(2p^{a-1})))^{(r)})$.
Now the $E_2$ term is a differential graded algebra in characteristic $p$,
so the $p$th power of an element of the even part passes to the next term
in the spectral sequence.
Therefore $E_3$ is also noetherian
 over its finitely
generated even degree part. As the spectral sequence stops, we get by
repeating this argument that $E_\infty$ is finitely generated.
But then so is the abutment.\qed

\begin{Problem}\label{integral}
There is a spectral sequence
$$E_1^{ij}=H^{i+j}(G,\gr_{-i}A)\Rightarrow H^{i+j}(G, A),$$
and the theorem says that $E_1$ is
finitely generated. The generators are in bidegree $(0,0)$
or $(i,j)$ with $i+j>0$, $i<0$. We would like this spectral
sequence to stop too.
It would suffice to know that if $J$ is a $G$-stable ideal in a ring like
$A$, then the
even part of $H^*(G,A/J)$ is integral over the image of the even part
of $H^*(G,A)$.
% \par View the spectral sequence as a module over its totally degenerate
% counterpart for the ring $\A$, where $\A$ describes the contracting
% family. It is what Grosshans calls $D$. You give $\A$ a rather stupid
% filtration, whose associated graded is $\A$ itself.
% Our interesting $E_1$ becomes a noetherian module over $\A$.\par
Of course this integrality would be implied by
finite generation of $H^*(G,S_A^*(A/J))$.
But recall  that  Nagata proves first that
$H^0(G,A/J)$ is integral over $H^0(G,A)$. So one could hope for a direct proof.
\end{Problem}

\begin{Problem}\label{divided}
If $M$ is a vector space, let $\Gamma^*(M)=\bigoplus_m(S^m(M^\#))^\#$
denote its divided power algebra.
It is probably too much to ask for a divided power structure on
the even part of the bigraded algebra
$H^*({\GL_n},\Gamma^*((\gl_n)(2)^{(1)}))$, extending the divided power
structure \cite[Appendix A2]{Eisenbud}, \cite{Roby} on the graded algebra
$H^0({\GL_n},\Gamma^*((\gl_n)(2)^{(1)}))$.
But such a structure would be
very helpful, as it would explain and enrich the supply of universal
cohomology
classes from \cite{Friedlander-Suslin}. And
$H^{i+j}(G,\gr_{-i}A)\Rightarrow H^{i+j}(G, A)$ would then undoubtedly
 stop. See our treatment of the rank one case below.
\end{Problem}

 \begin{Problem}
 To improve on the conditions of theorem \ref{fingrgr}, rather than on its
 conclusion, one should try and prove the following.
 If our finitely generated
 $A$ has a good filtration
 and $M$ is  a finite $A$-module on which $G$ acts compatibly,
 then $H^*(G,M)$ is finite  over $A^G$.
Of course this would again be implied by finite generation of
$H^*(G,S_A^*(M))$.
% \par An instructive case should be $A=\bigoplus_{n\geq0}\nabla(n\varpi_i)$.
 \end{Problem}

\begin{Remark}\label{other group}
Let $G$ be a semi-simple group defined over $k$, and let $V$
be a tilting $G$-module \cite{Donkin tilting} of dimension $n$.
Choose a basis in $V$ and assume that the representation is faithful,
so that $G$ can be identified with
a subgroup of $\SL_n$.
% \par As the representation has a good filtration the map $G\to \SL_n$
% can at most be inseparable in the center of $G$. As we can take
% fixed points under the center anyway, we do not mind.\par
Assume $p> n/2$.
% Note that $\wedge^iV=\wedge^{n-i}V^\#$.
Then $S^*((\gl_n)^\#)$ has a good filtration as a $G$-module
 \cite[4.3]{Andersen-Jantzen}.
Also assume for every fundamental weight $\varpi_i$ of
$G$ that the symmetric algebra
$S^*(\nabla(\varpi_i))$ has a good filtration.
This happens for instance when $p$ exceeds the
dimensions of the fundamental representations of $G$.
Let $A$ be a finitely generated commutative $k$-algebra on which $G$ acts
rationally by $k$-algebra automorphisms.
Then again $H^*(G,\gr A)$ is finitely generated as a
$k$-algebra.
The proof is the same as the proof of theorem \ref{fingrgr}.
\end{Remark}

\section{Divided powers}
Let $W_2(k)=W(k)/p^2W(k)$ be the ring of Witt vectors of length two over $k$,
see \cite[II \S6]{Serre}.
One has an extension of algebraic groups
$$1\to \gl_n^{(1)} \to \GL_n (W_2(k)) \to \GL_n(k)\to 1,$$
whence a cocycle class $e_1\in H^2(\GL_n,\gl_n^{(1)})$.
We call it the Witt vector class for $\GL_n$.
Analogously there are Witt vector extensions
and Witt vector classes for other groups that are originally
defined over the integers, and for Frobenius kernels in them.

\begin{Remark}\label{restrict to ga1}
Observe that if $x_\alpha:\Ga\to \SL_n$ is a root homomorphism
corresponding to the root $\alpha$, then the restriction of
$e_1$ to $\Ga$ is non-trivial. In fact one may restrict further
to the Frobenius kernel ${\Ga}_{1}$. This ${\Ga}_{1}$ acts trivially
on $\gl_n^{(1)}$, and the $T$-equivariant projection of
$\gl_n^{(1)}$ onto its $p\alpha$ weight space is thus also
${\Ga}_{1}$-equivariant. Altogether we get an image $\beta$ of $e_1$ in
$H^2(\Ga_{1},(\gl_n^{(1)})_{p\alpha})^T$. It is just the Witt vector class
of ${\Ga}_{1}$. It is well known that this class is
non-trivial, compare \cite[I 4.22, I 4.25]{Jantzen},
\cite[\S 6]{Friedlander-Suslin}. Also note that the image of $\beta$ in
$H^2(\Ga_{1},\gl_n^{(1)})$ under the map induced by
the inclusion of
$(\gl_n^{(1)})_{p\alpha}$ into $\gl_n^{(1)}$ is the same as the restriction of
$e_1$ to $\Ga_{1}$.
\end{Remark}

\begin{Lemma}
The Witt vector class for $\GL_n$ coincides with the universal cohomology
 class $e_1$
of
Friedlander and Suslin, up to a non-zero scalar factor.
\end{Lemma}

\paragraph{Proof}Using
 \cite[Remark 1.2.1, Corollary 3.13]{Friedlander-Suslin}
we see it suffices to take $r=j=1$ and $q=2$ in
\cite[Theorem 4.5]{Friedlander-Suslin}.
 \qed

\subsection{Divided powers in rank one.}
If $R$ is a commutative $k$-algebra, and $M$ is a
finite dimensional vector space over $k$,
then the divided power algebra
$R\otimes \Gamma^*(M)=\Gamma^*_R(R\otimes M)$ is a ring with divided powers
(\cite[Appendix A2]{Eisenbud}).
We write $\Delta_{i,j}$ for the component
$\Gamma^{i+j}(M)\to \Gamma^{i}(M)\otimes\Gamma^{j}(M)$
of the comultiplication
map $\Gamma^*(M)\to \Gamma^*(M)\otimes \Gamma^*(M)$.
So $\Delta_{i,j}$ is the dual of the multiplication map
$S^i(M^\#)\otimes S^j(M^\#)\to S^{i+j}(M^\#)$.
If $v\in M$ has divided powers $v^{[i]}\in \Gamma^{i}(M)$,
then $\Delta_{i,j}(v^{[i+j]})=v^{[i]}\otimes v^{[j]}$.
Put $G=\GL_n$ and define $T$, $B^+$, $B^-$ as usual.
Actually $n$ will be two, but we optimistically keep writing $n$.
Let $r\geq1$.
As $G_{r}$ acts trivially on $\gl_n^{(r)}$, there is a
divided power structure on $H^{\even}(G_{r},\Gamma^*(\gl_n^{(r)}))=
H^\even(G_{r},k)\otimes \Gamma^*(\gl_n^{(r)})$,
with the $m$-th divided power operation going from
$H^{2a}(G_{r},\Gamma^b(\gl_n^{(r)}))$ to
$H^{2am}(G_{r},\Gamma^{bm}(\gl_n^{(r)}))$
for $b\geq1$.

\paragraph{}
It would be nice to extend the next theorem to arbitrary $n$.
We do not know how to put a divided power structure on
$\bigoplus_m H^{2m}(\GL_n,\Gamma^{m}(\gl_n^{(1)}))$.
Nevertheless we feel the theorem is named appropriately.

\begin{Theorem}[Divided powers in rank one]\label{divided powers}
Let $n=2$. There are classes $c[m]\in H^{2m}(\GL_n,\Gamma^{m}(\gl_n^{(1)}))$
so that
\begin{enumerate}
\item $c[1]$ is the Witt vector class $e_1$,
\item $(\Delta_{i,j})_*(c[i+j])=c[i]\cup c[j]$ for $i,j\geq1$.
\end{enumerate}
\end{Theorem}

\paragraph{Proof}
Let $\alpha$ be the negative root, and let $x_\alpha :\Ga\to \GL_n$ be
its  root homomorphism.
By Kempf vanishing $H^{2m}(\GL_n,\Gamma^{m}(\gl_n^{(1)}))=
H^{2m}(\Ga,\Gamma^{m}(\gl_n^{(1)}))^T$, see \cite[II 4.7c]{Jantzen}.
So we restrict to $\Ga$ along $x_\alpha$.
As $\Ga$  acts trivially on
$\Gamma^*((\gl_n^{(1)})_{p\alpha})$, we also have a divided power
structure on $H^{\even}(\Ga,\Gamma^*((\gl_n^{(1)})_{p\alpha}))=
H^{\even}(\Ga,k)\otimes \Gamma^*((\gl_n^{(1)})_{p\alpha})$.
Take the $m$-th divided power  in
$H^{\even}(\Ga,\Gamma^*((\gl_n^{(1)})_{p\alpha}))^T$
of the Witt vector class of $\Ga$
 and map it to
$H^{2m}(\Ga,\Gamma^m(\gl_n^{(1)}))^T$ using the inclusion of
$(\gl_n^{(1)})_{p\alpha}$ into $\gl_n^{(1)}$. One lands in
$H^{2m}(B^-,\Gamma^m(\gl_n^{(1)}))\cong H^{2m}(G,\Gamma^m(\gl_n^{(1)}))$.
This gives the desired class $c[m]$.
Indeed $(\Delta_{i,j})_*(c[i+j])=c[i]\cup c[j]$ holds in the context of
the divided power algebra $H^{\even}(\Ga,\Gamma^*((\gl_n^{(1)})_{p\alpha}))=
H^{\even}(\Ga,k)\otimes\Gamma^*((\gl_n^{(1)})_{p\alpha})$
and thus the result follows by functoriality of $\Delta_{i,j}$
and of cup product.
% \par Some ramblings. First of all, the restriction of $c[m]$
% to $H^{2m}(\Ga_{1},\Gamma^m(\gl_n^{(1)}))^T$
% is the $m$-th divided power of the Witt vector class of $\Ga_{1}$.
% But for any $G$-module $M$ we have $H^{i}(G_{1},M^{(1)})^{(-1)}=
% \ind_{B}^G(H^{i}(B_1,k)^{(-1)}\otimes M)$ by
% \cite{Andersen-Jantzen}, or \cite[I 4.27, II 12.7]{Jantzen}.
% More specifically, if $I^*$ is a resolution of $M^{(1)}$ by injective
% $G$-modules,
% then $H^{*}(G_{1},M^{(1)})$ is the homology of $(I^*)^{G_{1}}$
% and $(I^*)^{G_{1}}=\ind_{B}^G(((I^*)^{B_1})^{(-1)})$.
% Both $((I^*)^{B_1})^{(-1)}$ and its homology are
% $\ind_{B}^G$-acyclic.
% Therefore the adjunction map
% $\ind_{B}^G(H^{i}(B_1,k)^{(-1)}\otimes M)\to
% H^{i}(B_1,k)^{(-1)}\otimes M$ agrees with the restriction map
% $H^{i}(G_{1},M^{(1)})\to H^{i}(B_1,M^{(1)})$.
% \par Consider the inclusions $k[G]^{G_1}\subset k[G]^{T_1\Ga_1}
% \subset k[G]$ and apply $\ind_{B}^G$. Then the adjunction
% $\ind_{B}^G(k[G]^{T_1\Ga_1})\to k[G]^{T_1\Ga_1}$ has as image a
% $G$-submodule of $k[G]^{T_1\Ga_1}$, which must be $k[G]^{G_1}$.
% So this adjunction map agrees with the inclusion map.\par
% It induces an
% isomorphism between modules of
% invariants $H^{i}(G_{1},\Gamma^m(\gl_n^{(1)}))^G\to
% H^{i}(B_1,\Gamma^m(\gl_n^{(1)}))^{B}$.
% Thus, to see that the restriction of $c[m]$ to
% $H^{2m}(G_{1},\Gamma^m(\gl_n^{(1)}))$ is the desired $G$-invariant,
% it suffices to observe that the
% restriction to $H^{2m}(B_1,\Gamma^m(\gl_n^{(1)}))=
% H^{2m}(\Ga_{1},\Gamma^m(\gl_n^{(1)}))^{T_1}$ is the desired
% $B$-invariant.
\qed

\subsection{Other universal classes}
If $M$ is a finite dimensional vector space over $k$ and
$r\geq1$, we have a natural
homomorphism between symmetric algebras $S^*(M^{\#(r)})\to S^{*}(M^{\#(1)})$
induced by the
map $M^{\#(r)}\to S^{p^{r-1}}(M^{\#(1)})$ which raises an element to the power
$p^{r-1}$. It is a map of bialgebras.
Dually we have the bialgebra map
$\pi^{r-1}:\Gamma^{p^{r-1}*}(M(1))\to \Gamma^{*}(M(r))$ whose kernel is the
ideal generated by $\Gamma^1(M^{(1)})$ through $\Gamma^{p^{r-1}-1}(M^{(1)})$.
So $\pi^{r-1}$ maps $\Gamma^{p^{r-1}a}(M(1))$ onto $\Gamma^{a}(M(r))$.

\begin{Notation} We now introduce analogues of the classes $e_r$ and
$e_r^{(j)}$ of Friedlander and Suslin \cite[Theorem 1.2,
Remark 1.2.2]{Friedlander-Suslin}.
We write $\pi^{r-1}_*(c[ap^{r-1}])\in
H^{2ap^{r-1}}(G,\Gamma^{a}(\gl_n^{(r)}))$
as $c_r[a]$. Next we get
$c_r[a]^{(j)}\in H^{2ap^{r-1}}(G,\Gamma^{a}(\gl_n^{(r+j)}))$
by Frobenius twist.
\end{Notation}

\begin{Lemma}\label{cup}
The $c_i[a]^{(r-i)}$ enjoy the following properties ($r\geq i\geq1$)
\begin{enumerate}
\item \label{cuphom} There is a homomorphism of algebras
$S^*(\gl_n^{\#(r)})\to H^{2p^{i-1}*}(G_{r},k)$
given on $S^a(\gl_n^{\#(r)})=H^0(G_{r},S^a(\gl_n^{\#(r)}))$
by cup product with the restriction of
$c_i[a]^{(r-i)}$.
\item \label{new e_r} For $r\geq1$ the restriction of
$c_r[1]$ to $ H^{2p^{r-1}}(G_{1},\gl_n^{(r)})$ is nontrivial,
so that $c_r[1]$ may serve as the universal class $e_r$ in
\cite[Thm 1.2]{Friedlander-Suslin}.
\end{enumerate}
\end{Lemma}

\paragraph{Proof}
\ref{cuphom}.
By theorem \ref{divided powers} we have
$(\Delta_{a,b})_*(c_{i}[a+b])=(\Delta_{a,b}\pi^{r-1})_*(c[(a+b)p^{r-1}])=
(\pi^{r-1}\otimes\pi^{r-1})_*(\Delta_{ap^{r-1},bp^{r-1}})_*(c[(a+b)p^{r-1}])=
(\pi^{r-1}\otimes\pi^{r-1})_*(c[ap^{r-1}]\cup c[bp^{r-1}])=
c_{i}[a]\cup c_{i}[b]$ and thus
$(\Delta_{a,b})_*(c_{i}[a+b]^{(r-i)})=
c_{i}[a]^{(r-i)}\cup c_{i}[b]^{(r-i)}$
 by pull back along a Frobenius homomorphism.
Put $R=H^\even(G_{r},k)$ and restrict from $G$ to $G_r$.
We view $H^\even(G_r,\Gamma^{*}(\gl_n^{(r)}))$ as
$\Gamma^*_R(R\otimes\gl_n^{(r)}))$.
Now the cup product agrees with the usual $R$-valued pairing between
$S^*(\gl_n^{\#(r)})$ and $\Gamma^*_R(R\otimes\gl_n^{(r)})$.
Thus if $x\in S^a(\gl_n^{\#(r)})$,
$y\in S^b(\gl_n^{\#(r)})$, then $(xy)\cup (c_{i}[a+b]^{(r-i)})=
(x\otimes y)\cup((\Delta_{a,b})_*(c_{i}[a+b]^{(r-i)}))=
(x\cup c_{i}[a]^{(r-i)})(y\cup c_{i}[b]^{(r-i)})$.
\\
\ref{new e_r}. In fact if we restrict $c_r[1]$ as in remark
\ref{restrict to ga1} to
$H^{2p^{r-1}}(\Ga_{1},(\gl_n^{(r)})_{p^r\alpha})=
H^{2p^{r-1}}(\Ga_{1},k)\otimes(\gl_n^{(r)})_{p^r\alpha}$, then even
that restriction is nontrivial. That is because the Witt
vector class generates the polynomial ring $H^\even(\Ga_{1},k)$, see
\cite[I 4.26]{Jantzen}.
\qed

\begin{Corollary}\label{dominate}Let $n=2$, $r\geq1$. Further
let $A$ be a
commutative $k$-algebra with $\SL_n$ action and $J$ an invariant ideal in $A$
so that the algebra $A/J$ is finitely generated.
Then $H^0(\SL_n,H^*((\SL_n)_r,A/J))$ is a noetherian module over a
finitely generated
subalgebra of $H^\even(\SL_n,A)$.
\end{Corollary}

\paragraph{Proof} We may assume $A$ is finitely generated.
Put $C=H^0((\SL_n)_r,A))$. Then $C$ contains the elements of $A$ raised to
the power $p^r$, so  $C$
is also a finitely generated algebra and $A/J$ is
a noetherian module over it. By \cite[Thm 1.5]{Friedlander-Suslin}
$H^*((\SL_n)_r,A/J)$ is a noetherian module over the finitely generated
algebra
$$R=\bigotimes_{a=1}^rS^*((\gl_n^{(r)})^\#(2p^{a-1}))\otimes
C.$$
Then, by invariant theory \cite[Thm. 16.9]{Grosshans book},
$H^0(\SL_n,H^*((\SL_n)_r,A/J))$ is a noetherian module over the finitely
generated algebra $H^0(\SL_n,R)$.
By lemma \ref{cup}
we may take the algebra homomorphism $R\to H^*((\SL_n)_r,A/J)$ of
\cite{Friedlander-Suslin}
to be based on cup product with our $c_i[a]^{(r-i)}$.
But then the map
$H^0(\SL_n,R)\to H^*((\SL_n)_r,A/J)$
factors, as  a linear map, through $H^\even(\SL_n,A)$.
\qed

\subsection{The contraction}
Let $A$ be a finitely generated commutative $k$-algebra on which $\SL_n$
acts rationally. Recall that $A$ comes with an increasing
filtration $A_{\leq i}=\sum_{\hgt(\lambda)\leq i}A_{\leq\lambda}$
whose associated graded is the ring $\gr A$. Let $\A$ be
the subring of the polynomial ring $A[t]$ generated by the subsets
$t^iA_{\leq i}.$
This ring $\A$, denoted $D$ by Grosshans,
is the coordinate ring of a flat family with general fibre $A$
and special fibre $\gr A$ \cite[Theorem 13]{Grosshans contr}.

\subsection{The special fibre}\label{special fibre}
There is a homomorphism of graded algebras
$\A\to \gr A$ with kernel $t\A$, mapping
$t^i\sum_{\hgt(\lambda)= i}A_{\leq\lambda}$ onto
$\sum_{\hgt(\lambda)=i}A_{\leq\lambda}/A_{<\lambda}$
by `dropping $t^i$'.
By corollary \ref{dominate} the even part of $E_2^{0*}$ in proposition
\ref{stops} is noetherian over a finitely generated subalgebra $R$ of
$H^\even(\SL_n,\A)$. Therefore proposition
\ref{stops} implies in the usual way
(\cite[Lemma 1.6]{Friedlander-Suslin})
that in fact
the abutment $H^*(\SL_n,\gr A)$ is noetherian over $R$.

\subsection{The general fibre}
One gets a homomorphism $\A\to A$  by substituting a nonzero scalar for $t$.
Let us use the substitution $t\mapsto 1$.
On $\A$ we put the filtration with
$\A_{\leq m}=\sum_{i\leq m}t^iA_{\leq i}$.
Then the associated graded $\gr\A$ is just $\A$ itself.
The map $\A\to A$ is compatible with the filtrations, so we get a map
of spectral sequences from
$$E(\A):\quad E_1^{ij}(\A)=H^{i+j}(G,\gr_{-i}\A)\Rightarrow H^{i+j}(G, \A)$$
to
$$E(A):\quad E_1^{ij}(A)=H^{i+j}(G,\gr_{-i}A)\Rightarrow H^{i+j}(G, A).$$
Note that by the construction of the Grosshans filtration
$H^{0}(G,\gr_{-i}A)$ vanishes for $i\neq0$.
Further $E_1^{**}(A)$
is finitely generated by theorem \ref{fingrgr}
and therefore there are for each $m$
only finitely many non-zero $E_1^{m-i,i}(A)$. This makes that,
even though the spectral sequence is a second quadrant
spectral sequence,
the abutment will be finitely generated as soon as $E_\infty^{**}(A)$ is.

\begin{Theorem}[Cohomological invariant theory in rank one]
Let $A$ be a finitely generated
commutative $k$-algebra on which $\SL_2$ acts rationally by algebra
automorphisms. Then $H^*(\SL_2,A)$ is finitely generated as a
$k$-algebra.
\end{Theorem}

\paragraph{Proof}
We combine the above.
The spectral sequence
$$E(\A):\quad E_1^{ij}(\A)=H^{i+j}(G,\gr_{-i}\A)\Rightarrow H^{i+j}(G, \A)$$
is pleasantly boring: It does not just degenerate,
even its abutment is the same as its $E_1$.
The spectral sequence $E(A)$ is a module over it \cite{Massey}.
% \par In the treatment of multiplicative structure for spectral sequences
%  of filtered complexes in \cite{Benson II}
%  one must allow the filtration to go
%  on negatively. Benson takes a descending filtration and
%  starts with the full complex. We do not start there.
%  To make it descend, we use negative indices, just like in the exercises
%  at the end of chapter XV in Cartan-Eilenberg.\par
 In particular,
$E(A)$ is a module over the finitely generated subring $R$ of
$H^\even(\SL_n,\A)$ over which $E_1^{**}(A)$ is noetherian by
\ref{special fibre}. (Yes, the
homomorphism $R\to E_1^{**}(A)$ agrees with the homomorphism
$E_1^{**}(\A)\to E_1^{**}(A)$.)
So the usual argument (see proof of \cite[Lemma 1.6]{Friedlander-Suslin} or
\cite[Lemma 7.4.4]{Evens book})
shows that $E(A)$ stops and that $E_\infty^{**}(A)$ is noetherian over $R$.
\qed

\begin{Corollary}
Let $A$ be a finitely generated
commutative $k$-algebra on which $\SL_2$ acts rationally by algebra
automorphisms. Assume that $A$ has a good filtration.
Let $M$ be a finitely generated $A$ module with a compatible $\SL_2$ action.
Then $M$ has finite good filtration dimension.
\end{Corollary}

\paragraph{Proof} Write $G=\SL_2$.
If we tensor $A$ with the multicone $k[G/U]$ of example \ref{multicone},
then the result is a finitely generated $G$-acyclic $k$-algebra
$A\otimes k[G/U]$ over which
we have a finitely generated module $M\otimes k[G/U]$.
As $H^*(G,S^*_{A\otimes k[G/U]}(M\otimes k[G/U]))$ is finitely generated
by the theorem,
$H^i(G,M\otimes k[G/U])$ vanishes for $i>>0$. But this means that
$M$ has finite good filtration dimension.
\qed

\end{document}